\documentclass[11pt]{amsart}

\usepackage{amssymb} 
\usepackage[mathscr]{eucal} 
\usepackage[all]{xy}
\usepackage{amsmath} 
\usepackage{epsfig}
\usepackage{amscd}

\newtheorem{TEO}{Theorem}[section]
\newtheorem{PROP}[TEO]{Proposition}
\newtheorem{LEM}[TEO]{Lemma}

\newtheorem{COR}[TEO]{Corollary}

\newcommand\dual{\mathrel{\raise3pt\hbox{$\underline{\mathrm{\thinspace d
\thinspace}}$}}}

\newcommand\proj{\mathbb P}
\newcommand\Z{\mathbb Z}

\def\Z{{\mathbb Z}}

\def\Z{{\mathbb Z}}

\begin{document}


\footnote{
Partially supported by PRIN 2007 MIUR: "Moduli, strutture geometriche e loro applicazioni " and PRIN 2006 of MIUR "Geometry on algebraic varieties".

AMS Subject classification: 14H10, 14J28. }

\title{On the second Gaussian map for curves on a K3 surface }

\author[E. Colombo]{Elisabetta Colombo}
\address{Dipartimento di Matematica,
Universit\`a di Milano, via Saldini 50,
     I-20133, Milano, Italy } \email{{\tt
elisabetta.colombo@unimi.it}, Fax:+39-02-50316090}

\author[P. Frediani]{Paola Frediani}
\address{ Dipartimento di Matematica, Universit\`a di Pavia,
via Ferrata 1, I-27100 Pavia, Italy } \email{{\tt
paola.frediani@unipv.it}, Fax:+39-0382-985602}

\maketitle

\begin{abstract}
By a theorem of Wahl, for canonically embedded curves which are
hyperplane sections of K3 surfaces, the first Gaussian map is not
surjective. In this paper we prove that if $C$ is a general
hyperplane section of high genus (greater than 280) of a general
polarized K3 surface, then the second Gaussian map of $C$ is
surjective. The resulting bound for the genus $g$ of a general
curve with surjective second Gaussian map is decreased to $g
>152$.
\end{abstract}

\section{Introduction}

The first Gaussian map, or Wahl map, for the canonical series has been extensively studied, and it has been shown that for a general curve of genus $\geq 10$, different from 11, it is surjective (\cite{chm}, \cite{voi}).

Wahl (\cite{wahl3}) has given a deformation theoretic
interpretation of the first Gaussian map, showing that if a
canonical curve can be extended in projective space as a
hyperplane section of a surface which is not a cone, then the
first Gaussian map is not surjective.

In particular in \cite{wahl3} it is proven that if a curve lies on a K3 surface, the first Gaussian map can't be surjective (see also \cite{bm}).

The obstruction to the surjectivity of the first Gaussian map for a curve in a K3 surface is given by the extension class of the cotangent sequence
 $$0\rightarrow K_C^{-1}\rightarrow {\Omega^1_X}_{|C}\rightarrow
K_C\rightarrow 0,$$
which is a non trivial element in the kernel of the dual of the first Gaussian map (see \cite{bm}).

This paper is concerned with the second Gaussian map, $\gamma^2_C:I_2(K_C)\rightarrow H^0(C,4K_C)$. In fact also the second Gaussian map encodes some interesting geometry. Our geometrical
motivation to study it comes from its relation with the curvature of the
moduli space  ${M_g}$ of curves of genus $g$ with the Siegel
metric induced by the period map $j:{ M_g}\rightarrow {A_g}$, that
we analyzed in \cite{cf2}.

There the curvature was computed using a formula for the
associated second fundamental form given in \cite{cpt}. In
particular in \cite{cpt} it is proven that the second fundamental
form lifts the second Gaussian map as stated in an
unpublished paper of Green and Griffiths (cf. \cite{green}).

In \cite{cf2}, (3.8) we gave a formula for the holomorphic
sectional curvature of ${M_g}$ along
 a Schiffer variation $\xi_P$, for $P$ a point  on the curve
 $C$, in terms of the holomorphic sectional curvature of ${A_g}$ and the second Gaussian map.

The relation of the second Gaussian map with curvature properties
of $M_g$ in $A_g$ suggested that its rank could give information on
the geometry of $M_g$ and of some sublocus of it.

Another interesting problem is to understand whether the surjectivity of second Gaussian map provides an obstruction to embed a curve in a surface as a hyperplane section.

In this paper we address this problem for curves in a K3 surface and we deduce results for the general curve.

Using cohomological techniques in the study of $\gamma^2_C$ for a curve in a K3 surface $X$ it is natural to consider the
"symmetric square" of the cotangent extension
 $$0\rightarrow {\Omega^1_X}_{|C}\otimes K_C^{-1}\rightarrow
S^2{\Omega^1_X}_{|C}\rightarrow K_C^2\rightarrow 0.$$
This does not give any obstruction to the surjectivity of $\gamma^2_C$ for the general curve in a general K3 surface, while it gives an obstruction if $C$ is any curve in an abelian surface (cf. \cite{cfp}). In fact in \cite{cfp} it is shown that if $C$ is a curve in an abelian surface $X$, then the corank of $\gamma^2_C$ is at least two.

In this paper (theorem (\ref{surj})) we prove surjectivity of the second Gaussian map for a general curve $C$  of high genus (for all $g > 280$) on a general polarized K3 surface.

This clearly implies surjectivity for the general curve of in the moduli space of curves
of any  genus $g >280$.

In corollary (\ref{general}) we decrease the lower bound for the genus of the general curve with surjective second Gaussian map till 152, using  examples given \cite{cf1}.

Note that, for dimensional reasons, surjectivity can be expected for a general curve of genus at least 18, therefore our bound is far from being optimal, nevertheless, it is the first known lower bound.

We recall that examples of curves whose second Gaussian map is
surjective were already given by in \cite{cf1} (for curves lying
on the product of two curves) and  in \cite{bafo}  (for complete
intersections). Note that using complete intersections it is not
possible to deduce surjectivity for the general curve of  any
sufficiently high genus, due to restrictions on the genus.
Moreover the first of such examples has much higher genus. The
examples of curves in a product of two curves are used in this
paper to decrease the bound.

On the other hand, theorem (\ref{surj}) shows that general curves on K3 surfaces of sufficiently high genus behave as  general curves in the moduli space, with respect to the second Gaussian map.

To prove our theorem we first show that sufficient conditions for the surjectivity of $\gamma^2_C$ for $C$ a curve in a K3 surface $X$ are given by the
surjectivity of the second Gaussian map
${\gamma^2_{{\mathcal O}_X(C)}}: I_2({\mathcal O}_X(C))\rightarrow H^0(S^2 \Omega^1_X \otimes {\mathcal O}_X(2C))$ and the vanishing of $H^1(S^2 \Omega^1_X \otimes {\mathcal O}_X(C))$.

To prove surjectivity of ${\gamma^2_{{\mathcal O}_X(C)}}$ we adapt the ideas used in \cite{clm3} to prove surjectivity of first Gaussian map ${\gamma^1_{{\mathcal O}_X(C)}}$.

More precisely, observe that a sufficient condition for the surjectivity of $\gamma^2_{{\mathcal O}_X(C)}$ is the vanishing of $H^1(  I^3_{\Delta_X} \otimes p^*({\mathcal O}_X(C)) \otimes q^*({\mathcal O}_X(C)))$, where $\Delta_X$ is the diagonal in $X \times X$ and $p,q$ are the two projections to $X$.

The idea is to consider the blow-up $Y$ of $X \times X$ along the diagonal $\Delta_X$ and to use Kawamata-Viehweg vanishing theorem (\cite{ka}, \cite{vi}) as follows. Let $E$ be the exceptional divisor and denote by $\pi: Y \rightarrow X \times X$ the natural morphism and by $f := p \circ \pi$, $g := q \circ \pi$. Then
$$H^1(  I^3_{\Delta_X} \otimes p^*({\mathcal O}_X(C)) \otimes q^*({\mathcal O}_X(C))) \cong H^1(Y,  f^*({\mathcal O}_X(C)) \otimes g^*({\mathcal O}_X(C))(-3E) )$$
$\cong H^1(Y,  f^*({\mathcal O}_X(C)) \otimes g^*({\mathcal O}_X(C)) \otimes K_Y (-4E)),$
since  $K_Y = {\mathcal O}_Y(E)$.\\

So by Kawamata-Viehweg vanishing theorem, it suffices to prove that
the line bundle $L:= f^*({\mathcal O}_X(C)) \otimes g^*({\mathcal O}_X(C))  (-4E)$ is big and nef. \\
Now notice that if one decomposes ${\mathcal O}_X(C)$ as $\otimes_{i=1}^4 A_i$, where $A_i$ are line bundles on $X$, then $L = \otimes_{i=1}^4 (f^*(A_i) \otimes g^*(A_i)(-E))$. To obtain that $L$ is big and nef, we ask suitable conditions on the line bundles $A_i$, and we study the sublinear system of $|f^*(A_i) \otimes g^*(A_i)(-E)|$ given by $\proj(\Lambda^2(H^0(A_i)))$ (cf. lemma (\ref{cilomi})).

The vanishing of $H^1(S^2 \Omega^1_X \otimes {\mathcal O}_X(C))$ relies on a similar argument, but it requires a more refined version of it.
In fact, given a decomposition of ${\mathcal O}_X(C)$ as ${\mathcal O}_X(D) \otimes {\mathcal O}_X(D')$, we have
$$H^1(S^2 \Omega^1_X \otimes {\mathcal O}_X(C)) =  H^1(X \times X, I^2_{\Delta_X}/I^3_{\Delta_X} \otimes p^*({\mathcal O}_X(D)) \otimes q^*({\mathcal O}_X(D'))),$$
hence its vanishing is implied by the vanishing of $H^1(X \times X, I^2_{\Delta_X} \otimes p^*({\mathcal O}_X(D)) \otimes q^*({\mathcal O}_X(D')))$ and of $H^2(X \times X, I^3_{\Delta_X} \otimes p^*({\mathcal O}_X(D)) \otimes q^*({\mathcal O}_X(D')))$.

So, with the same argument as above, it suffices to show that $f^*({\mathcal O}_X(D)) \otimes g^*({\mathcal O}_X(D'))  (-4E)$ is big and nef.
The strategy is now to choose ${\mathcal O}_X(D) = \otimes_{i=1}^4A_i$ and $D' =D + B$ with $B$ nef and effective, so we take $C \in |2D +B|$.

The above decompositions are shown on concrete examples of K3 surfaces $X$ and of curves $C$ in $X$, which are explicitely constructed via their Picard lattices (cf proposition \ref{picard}).

Finally, regardless the examples that we give, note that the conditions of the line bundles $A_i$ as in lemma (\ref{cilomi}) and  the decomposition ${\mathcal O}_X(C) =  {\mathcal O}_X(2D+B)$ force the genus of $C$ to be far from the optimal lower bound.

Anyhow, observe that the vanishing of $H^1(S^2 \Omega^1_X \otimes {\mathcal O}_X(C))$ itself, already implies that the curve $C$ must be of genus at least 31,  as one can check looking  at  the restriction of
$\Omega^1_X \otimes\Omega^1_X(C)$ to $C$ and the induced cohomology exact sequence.

\section{Preliminaries on Gaussian maps}

 Let $Y$ be a smooth complex projective variety and let $\Delta_Y\subset
Y\times Y$ be the diagonal. Let $L$ and $M$ be line bundles on $Y$.
For a non-negative integer $k$, the \emph{k-th Gaussian map}
associated to these data is the restriction to diagonal map
\begin{equation}\label{Gaussian1}\gamma^k_{L,M}:H^0(Y\times Y,I^k_{\Delta_Y}\otimes
L\boxtimes M )\rightarrow
H^0(Y,{I^k_{\Delta_Y}}_{|\Delta_Y}\otimes L\otimes M)\cong
H^0(Y,S^k\Omega_Y^1\otimes L\otimes M).
\end{equation}
Usually  \emph{first} Gaussian maps are simply referred to as
 \emph{Gaussian maps}. The exact sequence
\begin{equation}
\label{Ik} 0 \rightarrow I^{k+1}_{\Delta_Y} \rightarrow
I^k_{\Delta_Y} \rightarrow S^k\Omega^1_Y \rightarrow 0,
\end{equation}(where $S^k\Omega^1_Y$ is identified to its image via the diagonal map), twisted by $L\boxtimes M$, shows that the domain of the $k$-th
Gaussian map is the kernel of the previous one:
$$\gamma^k_{L,M}:
\ker \gamma^{k-1}_{L,M}\rightarrow H^0(S^k\Omega_Y^1\otimes L\otimes
M).$$ In our applications, we will exclusively deal with Gaussian
maps of order two, assuming also that the two line bundles
$L$ and $M$ coincide. For the reader's convenience, we spell out
these maps. The map $\gamma^0_L$ is the multiplication map of global
sections
\begin{equation}\label{symmetric}H^0(Y,L)\otimes
H^0(Y,L)\rightarrow H^0(Y,L^2)\end{equation}
 which obviously
vanishes identically on $\wedge^2 H^0(L)$.
 Consequently, $H^0(Y
\times Y, I_{\Delta_Y}\otimes L\boxtimes L)$ decomposes as
$\wedge^2 H^0(L)\oplus I_2(L)$, where $I_2(L)$ is the kernel of
$S^2H^0(Y,L)\rightarrow H^0(Y,L^2)$. Since $\gamma^1_L$ vanishes
on symmetric tensors, one writes
\begin{equation}\label{gamma1L}\gamma^1_L:\wedge^2H^0(L)\rightarrow H^0(\Omega^1_Y\otimes
L^2).\end{equation}
 Again, $H^0(Y\times Y, I^2_{\Delta_Y}\otimes L\boxtimes L)$
 decomposes as  the sum of $I_2(L)$ and the kernel of (\ref{gamma1L}). Since $\gamma_L^2$ vanishes identically on skew-symmetric
 tensors, one usually writes
 \begin{equation}\label{skew}\gamma^2_L:I_2(L)\rightarrow H^0(S^2\Omega_Y^1\otimes L^2)
 \end{equation}
  The primary object of this paper will be the
  second Gaussian map of the canonical line bundle
 $K_C$ on a curve $C$:

 $$\gamma^2_C:I_2(K_C)\rightarrow H^0(K_C^4)$$

In our situation $Y$ will be either a K3 surface $X$, or a smooth
irreducible projective curve $C$ on $X$.
\section{Main theorem}
\begin{TEO}
\label{surj}
If $X$ is a general polarized K3 surface of degree
$2g-2$ with $g>280$ and $C$ is a general hyperplane section of
$X$, then $\gamma^2_C$ is surjective.
\end{TEO}

Let us explain the strategy of the proof of theorem (\ref{surj}).
We have the following commutative diagram

\begin{equation}
\label{diagram}
\xymatrix{
I_2({\mathcal O}_X(C))\ar[dd]^r\ar[r]^{\gamma^2_{{\mathcal O}_X(C)} \ \ \ \ \ \ } & H^0(S^2 \Omega^1_X \otimes {\mathcal O}_X(2C)) \ar[dr]^{p_1} \\
&&H^0(S^2 \Omega^1_{X|C} \otimes K_C^2) \ar[dl]^{p_2}\\
I_2(K_C) \ar[r]_{\gamma^2_C}& H^0(K_C^4)}
\end{equation}
where $r$ and $p_1$ are restriction maps, and $p_2$ comes from the conormal extension.
More precisely, consider the exact sequence coming from the conormal extension
$$0 \rightarrow \Omega^1_{X|C} \otimes K_C \rightarrow S^2\Omega^1_{X|C} \otimes K^2_C \rightarrow K_C^4 \rightarrow 0,$$
then we have
$$H^0(S^2\Omega^1_{X|C} \otimes K^2_C) \stackrel{p_2}\rightarrow H^0(K_C^4) \rightarrow H^1(\Omega^1_{X|C} \otimes K_C) \cong H^0(T_{X|C})^*,$$
hence $p_2$ is surjective by the following lemma.
\begin{LEM}
\label{p2}
If $X$ is a general K3 surface and $C$ a general curve of genus at least 13 in the very ample linear system $|{\mathcal O}_X(C)|$ then $H^0(T_{X|C})=0$.
\end{LEM}
\proof
By the exact sequence given by restriction of $T_X$ to $C$, $H^0(T_{X|C})$ injects in $H^1(T_X(-C))$, which vanishes by lemma (2.3) of \cite{clm2}.
\qed\\

The theorem will follow if we prove that also the maps $\gamma^2_{{\mathcal O}_X(C)}$ and $p_1$ are surjective. In fact it suffices to exhibit examples of pairs $(X,C)$ where $X$ is a K3 and $C$ is a very ample curve in $X$ of any genus $g$ sufficiently high ($g \geq 281$) for which $\gamma^2_{{\mathcal O}_X(C)}$ and $p_1$ are surjective.
To do this we follow the strategy used in \cite{clm3} to study the first Wahl map.  More precisely, from the exact sequence
\begin{align}
\label{2/3}
&0 \rightarrow I^3_{\Delta_X} \otimes p^*({\mathcal O}_X(C)) \otimes q^*({\mathcal O}_X(C)) \rightarrow I^2_{\Delta_X} \otimes p^*({\mathcal O}_X(C)) \otimes q^*({\mathcal O}_X(C))  \\ \nonumber
& \rightarrow I^2_{\Delta_X}/I^3_{\Delta_X} \otimes p^*({\mathcal O}_X(C)) \otimes q^*({\mathcal O}_X(C)) \rightarrow 0
\end{align}
and taking global sections, we see that $\gamma^2_{{\mathcal O}_X(C)}$ is surjective if $H^1(  I^3_{\Delta_X} \otimes p^*({\mathcal O}_X(C)) \otimes q^*({\mathcal O}_X(C))) =0$.

The idea used in \cite{clm3} is to consider the blow-up $Y$ of $X \times X$ along the diagonal $\Delta_X$ and to use Kawamata-Viehweg vanishing theorem. Let $E$ be the exceptional divisor and denote by $\pi: Y \rightarrow X \times X$ the natural morphism and by $f := p \circ \pi$, $g := q \circ \pi$. Then
$$H^1(  I^3_{\Delta_X} \otimes p^*({\mathcal O}_X(C)) \otimes q^*({\mathcal O}_X(C))) \cong H^1(Y,  f^*({\mathcal O}_X(C)) \otimes g^*({\mathcal O}_X(C))(-3E) )$$
$\cong H^1(Y,  f^*({\mathcal O}_X(C)) \otimes g^*({\mathcal O}_X(C)) \otimes K_Y (-4E)),$
since  $K_Y = {\mathcal O}_Y(E)$.\\

So by Kawamata-Viehweg vanishing theorem, it suffices to prove that
$ f^*({\mathcal O}_X(C)) \otimes g^*({\mathcal O}_X(C))  (-4E)$ is big and nef. \\

Consider now the map
$$p_1: H^0(S^2 \Omega^1_X \otimes {\mathcal O}_X(2C))\rightarrow H^0(S^2 \Omega^1_{X|C} \otimes K_C^2).$$ Clearly $p_1$ is surjective if $H^1(S^2 \Omega^1_X \otimes {\mathcal O}_X(C)) =0$.

Our strategy to prove the surjectivity of $p_1$ is to adapt the
above idea for the vanishing of $H^1(Y,  f^*({\mathcal O}_X(C))
\otimes g^*({\mathcal O}_X(C)) \otimes (-3E))$ in order to show
that also $H^1(S^2 \Omega^1_X \otimes {\mathcal O}_X(C))$
vanishes.

To this end, let $H$ be a very ample divisor and assume that $C \in |2H + B|$, where $B$ is nef and effective. Then

\noindent $H^1(S^2 \Omega^1_X \otimes {\mathcal O}_X(C)) \cong H^1( X, q_*(I^2_{\Delta_X}/I^3_{\Delta_X} \otimes p^*({\mathcal O}_X(H)) \otimes q^*({\mathcal O}_X(H +B))))$ $ \cong H^1(X \times X, I^2_{\Delta_X}/I^3_{\Delta_X} \otimes p^*({\mathcal O}_X(H)) \otimes q^*({\mathcal O}_X(H +B)))$, where the last isomorphism comes from Leray spectral sequence.
So, by (\ref{2/3}), to prove surjectivity of $p_1$ it suffices to show that
$ H^1(  I^2_{\Delta_X} \otimes p^*({\mathcal O}_X(H)) \otimes q^*({\mathcal O}_X(H +B))) =0$ and
$ H^2(  I^3_{\Delta_X} \otimes p^*({\mathcal O}_X(H)) \otimes q^*({\mathcal O}_X(H +B))) =0.$ Using again the blow-up $Y$,  this is true if the line bundles $f^*({\mathcal O}_X(H)) \otimes g^*({\mathcal O}_X(H +B))(-4E)$ and $f^*({\mathcal O}_X(H)) \otimes g^*({\mathcal O}_X(H+B))(-3E)$ are big and nef.

In conclusion, if we prove that $f^*({\mathcal O}_X(H)) \otimes g^*({\mathcal O}_X(H +B))(-4E)$  and $f^*({\mathcal O}_X(H)) \otimes g^*({\mathcal O}_X(H +B))(-3E)$ are big and nef, then $p_1$ is surjective. Moreover also $ f^*({\mathcal O}_X(C)) \otimes g^*({\mathcal O}_X(C))  (-4E)$ is big and nef and therefore $\gamma^2_{{\mathcal O}_X(C)}$ is surjective.

Following \cite{clm3} we will exhibit pairs $(X,C)$ as above (where $C \in |2H +B|$) for which  $f^*({\mathcal O}_X(H)) \otimes g^*({\mathcal O}_X(H+B))(-4E)$ and  $f^*({\mathcal O}_X(H)) \otimes g^*({\mathcal O}_X(H+B))(-3E)$ are big and nef.
First of all observe that if there exist four line bundles $A_i$ $i=1,2,3,4$ on $X$ such that ${\mathcal O}_X(H) \cong A_1 \otimes A_2 \otimes A_3 \otimes A_4$, then we have
$$f^*({\mathcal O}_X(H)) \otimes g^*({\mathcal O}_X(H+B))(-4E) \cong$$
$$ \otimes_{i=1,2,3,4} (f^*(A_i) \otimes g^*(A_i)(-E)) \otimes g^*({\mathcal O}_X(B)),$$
and
$$f^*({\mathcal O}_X(H)) \otimes g^*({\mathcal O}_X(H+B))(-3E) \cong $$
$$\otimes_{i=1,2,3} (f^*(A_i) \otimes g^*(A_i)(-E)) \otimes (f^*(A_4) \otimes g^*(A_4 \otimes {\mathcal O}_X(B)).$$
These are big and nef under the conditions given in the following lemma.
\begin{LEM}
\label{cilomi}
Let $A_1$, $A_2$, $A_3$, $A_4$ be four base point free line bundles on a K3 surface $X$ with $A_j^2 \geq 2$, $j=1,2,3,4$ and such that $A_1$ is very ample. Assume that either
$A_2$, $A_3$, $A_4$ are very ample, or they define $(2:1)$ finite morphisms onto ${\proj}^2$ and that if $A_j^2 =2$, we have $(A_1\otimes A_2  \otimes A_3 \otimes A_4) \cdot A_j \geq 12$, then
$$\otimes_{i=1,2,3,4} (f^*(A_i) \otimes g^*(A_i)(-E))$$
is big and nef.

\end{LEM}
\proof
The proof is almost the same as the proof of lemma (2.2) of \cite{clm3}, but we reproduce it here for the reader's convenience.

If $A_i$ is very ample the linear system $|f^*(A_i) \otimes g^*(A_i)(-E)|$ on $Y$ has a sublinear system defining the morphism $F: Y \rightarrow Gr(1, \proj H^0(A_i)^*)$, associating to $(x,y) \in Y$ the line between $\phi_{A_i}(x)$ and $\phi_{A_i}(y)$, composed with the Pl\"ucker embedding. Notice that if $(x,y) \in E$, we can think of $(x,y)$ as a pair where $x \in X$, $y \in \proj T_{X,x}$, hence $F(x,y)$ is the line generated by $(d\phi_{A_i})_x(y)$. Therefore $f^*(A_i) \otimes g^*(A_i)(-E)$ is nef and it is also big, since the image of $X$ in $ \proj H^0(A_i)^*$ is not ruled.

Therefore, since $A_1$ is very ample, $(f^*(A_1) \otimes
g^*(A_1))(-E)$ is big, hence $\otimes_{i=1,2,3,4} (f^*(A_i)
\otimes g^*(A_i)(-E))$ is big and it can fail to be nef only on a
curve $Z$ contained in the indeterminacy locus of the maps $Y
\rightarrow Gr(1, \proj H^0(A_i)^*)$, $i =2,3,4$ if $A_i$ is not
very ample. Notice that $Z$ is a curve contained in $\{(x,y) \in
Y-E \ | \ \phi_{A_i}(x)= \phi_{A_i}(y) \} \cup \{(x,y) \in E \ | \
(d\phi_{A_i})_x(y) =0\}$. Assume that $A_i$ is not very ample,
hence, by assumption, it gives a $(2:1)$ morphism to  ${\proj}^2$.
If $Z \not\subset E$, let $\tau: X \times X \rightarrow X \times
X$ be the involution $\tau(x,y) =(y,x)$, then we can assume that
the image $\overline{Z}$ of $Z$ in $X \times X$ is such that
$\tau(\overline{Z}) = \overline{Z}$, because $\otimes_{i=1,2,3,4}
(f^*(A_i) \otimes g^*(A_i)(-E))$ is invariant under $\tau$. Then
the first (or second) projection $Z'$ of $\overline{Z}$ in $X$ is
$\phi_{A_i}^*(Z_1)$ for some curve $Z_1$ in ${\proj }^2$. If $L$
is a line in ${\proj}^2$ and $Z_1 \sim mL$ we have
$$\otimes_{i=1,2,3,4} (f^*(A_i) \otimes g^*(A_i)(-E))\cdot Z = 2(A_1 \otimes A_2 \otimes A_3 \otimes A_4) \cdot Z' - 4 E \cdot Z=$$
$$ 2(A_1 \otimes A_2 \otimes A_3 \otimes A_4) \cdot m \phi_{A_i}^*(L) - 4 E \cdot Z=2m(A_1 \otimes A_2 \otimes A_3 \otimes A_4) \cdot  A_i - 4 E \cdot Z,$$
therefore we are done if we show that $E \cdot Z = 6m$. Let $B$ be the ramification divisor of $\phi_{A_i}$, then $B$ is a smooth plane sextic and $E \cdot Z = mB \cdot L = 6m$, if the intersection of $E$ and $Z$ is transverse. This can be checked directly as in lemma (2.2) of \cite{clm3}.

If $Z \subset E$, then it is the strict transform of the ramification divisor $R$ on $X$ of $\phi_{A_i}$, hence $Z \cdot E = -c_1({\mathcal O}_{\proj T_X}(1)) \cdot Z = - degT_R = 18$. Therefore
$$\otimes_{i=1,2,3,4} (f^*(A_i) \otimes g^*(A_i)(-E)) \cdot Z =  6(A_1 \otimes A_2 \otimes A_3 \otimes A_4) \cdot A_i - 4 E \cdot Z \geq 0.$$
\qed

Let us now show the construction of the examples.

\begin{PROP}
\label{picard} There exist smooth $K3$ surfaces $X$ with Picard
lattice $\Gamma= \Z D \oplus \Z L \oplus \Z R \oplus \Z S \oplus
\Z T$ with intersection matrix $\bf{diag}$$(2h,-2k,-2j, -2l, -2m)$
with $j,k,l, m\geq 2$, $h \geq k+1, j+1, l+1, m+1$ and $D$ very
ample. Moreover $D+L$, $D+R$, $D+S$ and $D+T$ are base point free,
and they are either very ample, or they define $(2:1)$ morphisms
to ${\proj}^2$.
\end{PROP}
\proof
Observe that the lattice $\Gamma$ is even, nondegenerate and of signature $(1,4)$, hence it occurs as the Neron-Severi group of some algebraic K3 surface (cf. \cite{mo}, corollary 2.9).

We will show that there does not exist a class $F \in \Gamma$ such that $F^2 =-2$, $D \cdot F =0$. By well known results on periods of K3 surfaces (see e.g. \cite{beau}), this implies that there exists a $K3$ surface with Picard lattice $\Gamma$ and such that $D$ is ample.

Assume that $F = aD + bL +cR +dS + eT$ ($a,b,c,d, e \in \Z$) is such that $F \cdot D =0$, and $F^2 = -2$. The first equality implies $a =0$, and the second one yields $1 = k b^2 + j c^2 + ld^2 + m e^2$, which is absurd, since $k,j,l,m\geq 2$.

So $D$ is ample and $D^2 \geq 4$, hence $D$ is very ample provided that there does not exist an irreducible curve $F$ such that $F^2 = 0$, $F \cdot D =1,2$ (cf. \cite{sd}, or \cite{mori} theorem 5). But this cannot happen, because if we write $F = aD + b L + cR + dS + eT$, $a,b,c,d,e \in \Z$, $F \cdot D= 2ha \neq 1,2$ since $h \geq 2$. So $D$ is very ample.

Note that $(D+L)^2 = 2h -2k \geq 2$, since $h \geq k +1$.

First of all we show that for any $(-2)$-curve $F$, $F \cdot(D+L)
> 0$, hence $D+L$ is ample and it is base point free, provided
that there does not exist irreducible curves $F,G$ and an integer
$a \geq 2$ such that $D+L \sim aF +G$, with $F^2=0$, $G^2 =-2$, $F
\cdot G =1$ (cf. \cite{sd}, or \cite{mori} theorem 5). This
clearly cannot happen, since the product of two classes is always
even.

Set $F = aD + bL + cR +dS + eT$, $a, b,c,d,e \in \Z$ with $F^2 =-2$. Since $D$ is ample, $D \cdot F = 2h a >0$, hence $a >0$. $F^2 =-2$ yields $a^2 = \frac{k b^2  + j c^2 +ld^2+ me^2-1}{h}$.

 If $(D +L) (aD + bL + cR +dS+eT) = 2ah -2kb \leq 0$, we have $bk \geq ah>0$, so $b>0$ and
 $b^2 k^2 \geq h^2 a^2 = h  (b^2 k + c^2 j+d^2 l + e^2m-1)$. Thus we get
 $$b^2(k^2 -hk) -hc^2j -hd^2 l  -he^2m +h \geq 0,$$
 so if we set $h = k +1 +t$, $t\geq 0$, we obtain
 $$ 0 \leq b^2(k^2 -hk) -hc^2j -hd^2l -he^2m +h \leq b^2(k^2 -hk) +h =$$
 $$ t(1-b^2k) + k(1-b^2) +1 \leq -t +1$$
 hence we must have either $t =0$ or $t =1$.
 But if $t =0$, $h = k+1$, so
 $$0 \leq  b^2(k^2 -hk) -hc^2j -hd^2l -h e^2m+h =$$
 $$ k (1 -b^2) - (k+1)(c^2j +d^2l +e^2m) +1 \leq k(1-b^2) +1,$$
 so $b =1$, but then $ - (k+1)(c^2j+ d^2l  +e^2m) +1 \geq 0$, which is absurd.

 If $t =1$, $h =k+2$, so we have
 $$0 \leq -2b^2k + k+2 -(k+2)(c^2j + d^2l +e^2m) \leq k(1-2b^2) +2,$$
 thus we must have $b =1$, hence $0 \leq -k +2 -(k+2)(c^2j +d^2l+ e^2m) \leq -(k+2) (c^2j +d^2l+ e^2m)$, which implies $c =d=e=0$. But then $a^2 =   \frac{k  -1}{k+2}$, which is absurd.
So $D+L$ is ample and base point free. If $(D+L)^2 \geq 4$ a similar computation shows that there does not exists a curve $F$ such that $F^2 =0$ and $F \cdot (D+L) =1,2$, therefore $D+L$ is very ample.

If $(D+L)^2 =2$, it clearly defines a (2:1) morphism to ${\proj}^2$.

The same holds for $D+R$, $D+S$ and $D+T$.
\qed\\

{\em Proof} of (\ref{surj}).
Consider the K3 surfaces constructed in (\ref{picard}).

Let us set $A_1 = D$, $A_i$ $i = 2,3,4$ be either equal to $D$, or to $D+L$, or to $D+R$, or to $D+S$, or to $D+T$. Setting $H = A_1 + A_2 +A_3 +A_4 = 4D + aL +bR +cS+dT$, with $a,b,c,d \geq 0$, $a+b+c+d \leq 3$, $\tilde{H} = 2H +B$, where $B = nD + m(D+L) + r(D+R) +s(D+S)+t(D+T)$, with $m,n,r,s,t \geq 0$, lemma (\ref{cilomi}) applies. In fact it suffices to check that $H \cdot (D+L) \geq 12$, $H \cdot (D+R) \geq 12$, and $H \cdot (D+S) \geq 12$, $H \cdot (D+T) \geq 12$, which is true for our choices of $H$.  Hence $\otimes_{i=1,2,3,4} (f^*(A_i) \otimes g^*(A_i)(-E) \otimes g^*({\mathcal O}_X(B))$ is big and nef and  $(\otimes_{i=1,2,3} (f^*(A_i) \otimes g^*(A_i)(-E))) \otimes (f^*(A_4) \otimes g^*(A_4)\otimes ({\mathcal O}_X(B)) $ is also big and nef. If we take $C \in |\tilde{H}|$, we also have $H^1(S^2 \Omega^1_X \otimes {\mathcal O}_X(C)) =0$, thus $\gamma^2_C$ is surjective.

Now we have to check that with our choices of a non divisible
$\tilde{H}$, we obtain all the genera $g(C)=  1 + \frac{1}{2}
{\tilde{H}}^2\geq 281$ for curves $C \in |\tilde{H}|$. To this
end, it suffices to take $\tilde{H} =aD + s L + tR +v S + r T,$
where $s,t,v,r \geq 0$, relatively prime, $s+t+v+r\leq a-2$,
$9\leq a\leq 14$ and, if $a=9$, at most one among $s,t,v,r$ is
odd, if $a=10$,
 at most two
among $s,t,v,r$ are odd, if $a=11$, at most three among $s,t,v,r$
are odd.

In particular, let us start with $\tilde{H} = 9D +6L+R,$  for $C \in |\tilde{H}|$, we have
$$g(C) = 1 + \frac{1}{2} {\tilde{H}}^2 = 1 + 81h -36k -j.$$
Let us set $n = k-2$, $m = j-2 $, $h = n+3 +t$, where $n,m,t \geq 0$, then we have $h \geq m+3$, so $t \geq m-n$, $g(C) = 1 + 81(n+3+t) -36(n+2) -(m+2),$ and we have two cases:
\begin{enumerate}
\item $\rho := n-m \geq 0$, $t \geq 0$, then
$$g(C) = 1 + 81(m + \rho +3+t) -36(m + \rho +2) -(m+2)=170 + 45 \rho + 44m +81t,$$
with $t,m,\rho \in \Z$, $t,m, \rho \geq 0$.
\item $\alpha:= m-n \geq 0$, $t = \alpha + \beta$, with $\beta \geq 0$, then
$$g(C) = 1 + 81(n + 3+\alpha + \beta) -36(n+2) -(n+2 + \alpha)=170 +  44n +80 \alpha + 81\beta,$$
with $n, \alpha, \beta \in \Z$, $n, \alpha, \beta \geq0$.

\end{enumerate}

Since 44 and 45 are relatively prime it is clear that with $g(C)
=170 + 45 \rho + 44m +81t$ one gets all sufficiently high genera.
Using (1) and (2),  one can simply check that $g(C)$ runs through
all the integers greater than 620 and with the other choices of
$\tilde{H}$ one gets all genera $g$ greater than 280 except for $g
=321$.

For $g =321$ we consider the $K3$ surface contructed in
\cite{clm3} proposition (3.2) with Picard lattice given by
$\Gamma= \Z D \oplus \Z L$ with  $D^2= 4$, $L^2=2$, $D\cdot L=7$.
In \cite{clm3} it is proven that $D$ is very ample and $L$ defines
a $2:1$ finite morphism onto $\proj^2$. So if we set $A_i =D$, $i
=1,2,3,4$, $H = 4D$, $B = 3D+L$, $\tilde{H} = 2H + B = 11D +L$,
since $H \cdot L = 4D \cdot L = 28$, lemma (\ref{cilomi}) applies.
Hence $\otimes_{i=1,2,3,4} (f^*(A_i) \otimes g^*(A_i)(-E) \otimes
g^*({\mathcal O}_X(B))$ is big and nef and  $(\otimes_{i=1,2,3}
(f^*(A_i) \otimes g^*(A_i)(-E))) \otimes (f^*(A_4) \otimes
g^*(A_4)\otimes ({\mathcal O}_X(B)) $ is also big and nef. So, as
above, if we take $C \in |\tilde{H}|$, we also have $H^1(S^2
\Omega^1_X \otimes {\mathcal O}_X(C)) =0$, thus $\gamma^2_C$ is
surjective. Now it suffices to check that $C$ has genus $g = 1 +
\frac{1}{2} {\tilde{H}}^2 = 321.$

\qed

\begin{COR}
\label{general}
For the general curve of genus greater than 152, the second Gaussian map $\gamma^2_C$ is surjective.
\end{COR}
\proof By theorem (\ref{surj}) and semicontinuity of the corank of
$\gamma^2_C$,  for a general curve of genus greater than 280
$\gamma^2_C$ is surjective. Surjectivity for the general curve of
genus $153\leq g \leq 280$ can be proved exhibiting examples of
curves of genus $g$ with surjective second Gaussian map, which are
either hyperplane sections of a polarized K3 surface as in the
proof of (\ref{surj}), or in the product of two curves as in
\cite{cf1} theorem 3.1.

More precisely let $C_1$, $C_2$ be two smooth curves of respective
genera $g_1$, $g_2$, choose  divisors $D_i$ on $C_i$ of degree
$d_i$, $i=1,2$. Set $X = C_1 \times C_2$, let $C \in |
{p_1}^*(D_1) \otimes {p_2}^*(D_2)|$ be a smooth curve, where $p_i$
is the projection from $C_1 \times C_2$ on $C_i$, then $g(C) = 1 +
(g_2-1)d_1 + (g_1-1)d_2 + d_1 d_2$.

In \cite{cf1} we proved that if either $g_1,g_2 \geq 2$, $d_i \geq 2g_i + 5$, $i=1,2$, or $g_1 \geq 2$, $g_2 =1$, $d_1 \geq 2g_1 + 5$, $d_2 \geq 7$, or  $g_2 =0$, $d_2 \geq 7$, $d_2(g_1 -1) > 2d_1 \geq 4g_1 + 10$, then  $\gamma^2_C$ is surjective.

Then one has to check directly that these values of $g(C)$ cover all the remaining integers between 153 and 280.
\qed

\end{document}